\documentclass[11pt,a4paper]{article}

\usepackage{amsmath,amssymb}
\newtheorem{theorem}{Theorem}[section]
\newtheorem{corollary}{Corollary}[section]
\newtheorem{lemma}{Lemma}[section]
\newtheorem{proposition}{Proposition}[section]
\newtheorem{definition}{Definition}[section]
\newtheorem{remark}{Remark}[section]
\numberwithin{equation}{section}
\newtheorem{example}{Example}[section]
\newenvironment{proof}[1][Proof]{\noindent\textbf{#1.} }{\hfill {$\Box$}}
\numberwithin{equation}{section}

\begin{document}

\title{\textsc{On Exponential Stability for Skew-Evolution Semiflows on Banach Spaces }}
\author{\textsc{Codru\c{t}a Stoica} \\
Institut de Math\' ematiques  \\
Universit\' e Bordeaux 1  \\
France  \\
e-mail: \texttt{codruta.stoica@math.u-bordeaux1.fr}
}
\date{}
\maketitle

{\footnotesize \noindent \textbf{Abstract.} The paper emphasizes
the property of stability for skew-evolution semiflows on Banach
spaces, defined by means of evolution semiflows and evolution
cocycles and which generalize the concept introduced by us in
\cite{MeStBu_CJM}. There are presented several general
characterizations of this asymptotic property out of which can be
deduced well known results of the stability theory. A unified
treatment in the uniform and in the nonuniform setting is
given. The main results are also formulated in discrete time.}

{\footnotesize \vspace{3mm} }


{\footnotesize \vspace{2mm} }


\section{Preliminaries}

In recent years, an impressive progress concerning the study of
asymptotic behaviors for evolution equations can be emphasized,
which led to a vast literature due mostly to asymptotic properties
of linear operators semigroups. The possibility of reducing the
nonautonomous case in the study of evolutionary families or
skew-product flows to the autonomous case of evolution semigroups
on various Banach function spaces is considered an important way
to interesting applications. The study of the asymptotic behavior
of linear skew-product semiflows has been used in the theory of
evolution equations in infinite dimensional spaces. The approach
from point of view of asymptotic properties for the evolution
semigroup associated to the linear skew-product semiflows was
essential.

Interesting results on the stability of solutions for non
autonomous differential equations in Banach spaces were obtained
by E.A. Barbashin in 1967 and presented in \cite{Ba_NAU}. One of
the most important results of the stability theory is due to R.
Datko who had proved in 1970 in \cite{Da_JMAA} that a strongly
continuous semigroup of operators $\mathcal{S}=\{S(t)\}_{t\geq 0}$
of $\mathcal{C}_{0}$-class defined on a complex Hilbert space $H$
is uniformly exponentially stable if and only if for each vector
$x\in H$ the mapping given by $t\rightarrow \left\Vert S(t)x
\right\Vert$ is in $\mathcal{L}^{2}(\mathbb{R}_{+})$.

Later, in 1983, A. Pazy generalizes the previous result in
\cite{Pa_SV} for $\mathcal{L}^{p}(\mathbb{R}_{+})$, $p\geq 1$, in
the case of $\mathcal{C}_{0}$-semigroups in Banach spaces.

Among the remarkable results concerning the property of stability
is the classic theorem of R. Datko who proves in 1972 in
\cite{Da_JMA} that a evolutionary process denoted
$\mathcal{U}=\{U(t,s)\}_{t\geq s \geq 0}$, with uniform
exponential growth, is uniformly exponentially stable if and only
if there exists an exponent $p\geq 1$ such that
\begin{equation}
\underset{s\geq 0}{\sup }\int_{s}^{\infty }\left\Vert U
(t,s)x\right\Vert ^{p}dt<\infty , \ \forall x\in  V,
\end{equation}
$ V$ being a Banach space.

In 1986, S. Rolewicz presents in \cite{Ro_JMAA} a similar
characterization for the uniform exponential stability property,
proving that if for a continuous nondecreasing mapping
$F:\mathbb{R}_{+}\rightarrow\mathbb{R}_{+}$ with the properties
$F(0)=0$ and $F(t)>0, \ \forall t>0$ and for an evolutionary
process $\mathcal{U}=\{U(t,s)\}_{t\geq s \geq 0}$ on a Banach
space $ V$, with exponential growth, following relation holds
\begin{equation}
\underset{s\geq 0}{\sup }\int_{s}^{\infty }F(\left\Vert U
(t,s)x\right\Vert) dt<\infty , \ \forall x\in  V,
\end{equation}
then the evolutionary process is uniformly exponentially stable.

It was J.M.A.M. van Neerven who, in 1995, referring to the case of
a $\mathcal{C}_{0}$-semigroup, proposes in \cite{Ne_JFA}
characterizations of the uniform exponential stability in terms of
functions spaces. It is proven that a semigroup denoted
$\mathcal{S}=\{S(t)\}_{t\geq 0}$ is exponentially stable if there
exists a Banach function space $E$ over $\mathbb{R}_{+}$ with the
property $\underset{t\rightarrow \infty}{\lim }\left\Vert
\chi_{[0,t]}\right\Vert_{E}=\infty$ such that $\left\Vert
S(\cdot)x\right\Vert \in E, \ \forall x\in V$, $ V$ being a Banach
space. The result concerns the evolution semigroup associated to
an evolution operator.

In \cite{Ne_IEOT} J.M.A.M. van Neerven extends results obtained
previously by A. Pazy, S. Rolewicz and R. Datko, as, for instance,
in the case of lower semicontinuous functionals. There are
presented implications of the fact that a
$\mathcal{C}_{0}$-semigroup is not uniformly exponentially stable.

In \cite{Bu_SED} C. Bu\c{s}e proves a sufficient condition for the
uniform exponential stability of a strongly continuous semigroup
on a Banach space, as simplified variant of a theorem of J.M.A.M.
van Neerven. A nonautonom version of the theorem for evolution
families is also presented and it is emphasized the fact that it
implies results obtained by R. Datko, A. Pazy and S. Rolewicz.
There can also be emphasized the results obtained by P. Preda, A.
Pogan and C. Preda in \cite{PrPoPr_CzMJ}, where it is
characterized, from continuous and discrete point of view, the
exponential stability of an evolution process in terms of
existence of functionals on certain function spaces.

A generalization of the result of J.M.A.M. van Neerven and
characterizations of Perron type for evolution operators were
given by C. Preda and P. Preda, who in \cite{PrPr_MIR07} proves a
necessary and sufficient condition for an evolution process of
class $\mathcal{C}(0,e)$ to be uniformly exponentially stable and
by M. Megan, A.L. Sasu and B. Sasu in \cite{MeSaSa_RMA} evolution
operators on Banach spaces, respectively in \cite{MeSaSa_MM} where
there are emphasized discrete and continuous characterizations for
the uniform exponential expansiveness of skew-product flows, by
means of uniform complete admissibility of the pairs
$(\textit{c}_{0}(\mathbb{N},  V), \textit{c}_{0}(\mathbb{N}, V))$
respectively $(\mathcal{C}_{0}(\mathbb{R}_{+},
V),\mathcal{C}_{0}(\mathbb{R}_{+}, V))$, $ V$ being a Banach
space.

In \cite{MeSaSa_BBMS}, M. Megan, A.L. Sasu and B. Sasu prove
necessary and sufficient conditions for the uniform exponential
stability of evolution equation on Banach spaces, by emphasizing
techniques from the domain of skew-product semiflows and of some
function spaces. There are obtained generalizations of some
results by R. Datko, S. Rolewicz and J.M.A.M. van Neerven.

In \cite{StMe_UVT}, we have studied, among other asymptotic
properties, the stability of evolution operators in the uniform
setting, by emphasizing relations between the exponential
stability and integral stability.

The concept of skew-evolution semiflow, defined by us by means of
evolution semiflows and evolution cocycles seems to be more
appropriate for the study of asymptotic behaviors of evolution
equations by means of evolution operators and also for the
nonuniform setting. In fact, we prove that the evolution operators
are a particular case of evolution cocycles and that the
skew-evolution semiflows are generalizations of the skew-product
semiflows.

\section{Definitions and examples}

Let $( X,d)$ be a metric space, $ V$ a Banach space, $ V^{*}$ its
topological dual, $\mathcal{B}( V)$ the space of all $ V$-valued
bounded operators defined on $ V$. We denote $ Y= X\times V$,
$ T=\left\{ (t,t_{0})\in \mathbb{R}%
_{+}^{2}:t\geq t_{0}\right\}$ and $\Delta=\left\{ (m,n)\in
\mathbb{N}^{2}:m\geq n\right\}$. The norm of vectors on $ V$ and
on $ V^{*}$ and of operators on $\mathcal{B}( V)$ is denoted by
$\left\Vert \cdot \right\Vert$. Let $I$ be the identity operator
on $ V$. We will consider the set
\[
\mathcal{R}=\{R:\mathbb{R}_{+}\rightarrow\mathbb{R}_{+}| \ R \
\textrm{nondecreasing}, \ R(0)=0, \ R(t)>0, \ \forall t>0\}.
\]

\begin{definition}\rm\label{d_semiflow}
A mapping $\varphi: T\times  X\rightarrow  X$ is called
\textit{evolution semiflow} on $ X$ if it satisfies the following
properties

$(es_{1})$ $\varphi(t,t,x)=x, \ \forall (t,x)\in
\mathbb{R}_{+}\times  X$

$(es_{2})$ $\varphi(t,s,\varphi(s,t_{0},x))=\varphi(t,t_{0},x), \
\forall (t,s),(s,t_{0})\in  T, \ \forall x\in
 X$.
\end{definition}

\begin{definition}\rm\label{d_cociclu_2}
A mapping $\Phi: T\times  X\rightarrow \mathcal{B}( V)$ that
satisfies the following properties

$(ec_{1})$ $\Phi(t,t,x)=I, \ \forall t\geq0,\ \forall x\in
 X$

$(ec_{2})$
$\Phi(t,s,\varphi(s,t_{0},x))\Phi(s,t_{0},x)=\Phi(t,t_{0},x),
\forall (t,s), (s,t_{0})\in  T,\forall x\in  X$

\noindent is called \emph{evolution cocycle} over the evolution
semiflow $\varphi$.
\end{definition}

\begin{definition}\rm\label{d_lses}
A mapping $C: T\times  Y\rightarrow  Y$ defined by
\begin {equation}
C(t,s,x,v)=(\varphi(t,s,x),\Phi(t,s,x)v), \ \forall (t,s,x,v)\in
 T\times  Y
\end{equation}
where $\Phi$ is an evolution cocycle over the evolution semiflow
$\varphi$, is called \emph{skew-evolution semiflow} on $ Y$.
\end{definition}

\begin{example}\rm\label{ex_cev}
Let $f:\mathbb{R}\rightarrow\mathbb{R}_{+}$ be a function which is
nondecreasing on the interval $(-\infty,0)$ and decreasing on the
interval $(0,\infty)$ with the property that there exists
\[
\underset{t\rightarrow\pm\infty}{\lim }f(t)=l\in(0,\infty).
\]
Let us consider the set
$\mathcal{C}=\mathcal{C}(\mathbb{R},\mathbb{R})$ of all continuous
functions given by $x:\mathbb{R}\rightarrow \mathbb{R}$, endowed
with the uniform convergence topology on compact subsets of
$\mathbb{R}$. $\mathcal{C}$ is metrizable by respect to the metric
\begin{equation*}
d(x,y)=\sum_{n=1}^{\infty}\frac{1}{2^{n}}\frac{d_{n}(x,y)}{1+d_{n}(x,y)},
\ \textrm{where} \ d_{n}(x,y)=\underset{t\in
[-n,n]}\sup{|x(t)-y(t)|}.
\end{equation*}

Let $ X$ be the closure of the set $\{{f_{t}, \ t\in
\mathbb{R}\}}$ in $\mathcal{C}$, where we consider
$f_{t}(\tau)=f(t+\tau), \ \forall \tau\in \mathbb{R}$. Then $(
X,d)$ is a metric space and the mapping
\begin{equation*}
\varphi: T\times  X\rightarrow  X, \ \varphi(t,s,x)=x_{t-s}
\end{equation*}
is an evolution semiflow on $ X$.

Let $ V=\mathbb{R}^{3}$ be a Banach space with the norm
\[
\left\Vert (v_{1},v_{2},v_{3})\right\Vert=|v_{1}|+|v_{2}|+|v_{3}|.
\]
The mapping $\Phi: T\times  X\rightarrow \mathcal{B}( V)$ given by
\begin{equation*}
\Phi(t,s,x)(v_{1},v_{2},v_{3})=\left(
e^{\alpha_{1}\int_{s}^{t}x(\tau-s)d\tau}v_{1},e^{-\alpha_{2}\int_{s}^{t}
x(\tau-s)d\tau}v_{2},e^{l+\alpha_{3}\int_{s}^{t}x(\tau-s)d\tau}v_{3}\right)
\end{equation*}
where $(\alpha_{1},\alpha_{2},\alpha_{3})\in\mathbb{R}^{3}$, is an
evolution cocycle on $\mathbb{R}^{3}$ and $C=(\varphi,\Phi)$ is a
skew-evolution semiflow on $ Y$.
\end{example}

Some particular classes of measurable skew-evolution semiflows
will be defined in what follows, useful for various
characterizations of the property of stability in the uniform as
well as in the nonuniform case.

\begin{definition}\rm\label{def_taremas}
A skew-evolution semiflow $C=(\varphi,\Phi)$ is called
\emph{strongly measurable} if for all
$(t_{0},x,v)\in\mathbb{R}_{+}\times Y$ the mapping
\begin{equation}
s\mapsto\left\Vert\Phi(s,t_{0},x)v\right\Vert
\end{equation}
is measurable on $[t_{0},\infty)$.
\end{definition}

\begin{definition}\rm\label{def_startaremas}
A skew-evolution semiflow $C=(\varphi,\Phi)$ is called
\emph{$*$-strongly measurable} if for all $(t,t_{0},x,v^{*})\in
T\times X\times V^{*}$ the mapping
\begin{equation}
s\mapsto\left\Vert\Phi(t,s,\varphi(s,t_{0},x))^{*}v^{*}\right\Vert
\end{equation}
is measurable on $[t_{0},t]$.
\end{definition}

\section{Uniform exponential stability for skew-evolution semiflows}

Some definitions of asymptotic properties in the uniform setting
are given in what follows.

\begin{definition}\rm\label{ueg}
A skew-evolution semiflow $C =(\varphi,\Phi)$ has \emph{uniform
exponential growth} if there exist some constants $M\geq 1$ and
$\omega>0$ such that following relation holds
\begin{equation}
\left\Vert \Phi(t,t_{0},x)v\right\Vert \leq Me^{\omega
(t-s)}\left\Vert \Phi(s,t_{0},x)v\right\Vert,
\end{equation}%
for all $(t,s),(s,t_{0})\in  T$ and all $(x,v)\in
 Y$.
\end{definition}

\begin{definition}\rm\label{us}
A skew-evolution semiflow $C =(\varphi,\Phi)$ is called
\emph{uniformly stable} if there exists a constant $N\geq 1$ such
that
\begin{equation}
\left\Vert \Phi(t,t_{0},x)v\right\Vert \leq N\left\Vert
\Phi(s,t_{0},x)v\right\Vert
\end{equation}%
\noindent for all $(t,s),(s,t_{0})\in  T$ and all $(x,v)\in  Y$.
\end{definition}

\begin{definition}\rm\label{ues}
A skew-evolution semiflow $C =(\varphi,\Phi)$ is said to
be \emph{uniformly exponentially stable} if there exist some constants $N\geq1$ and $%
\nu >0$ such that
\begin{equation}
\left\Vert \Phi(t,t_{0},x)v\right\Vert \leq Ne^{-\nu
(t-s)}\left\Vert \Phi(s,t_{0},x)v\right\Vert
\end{equation}%
\noindent for all $(t,s),(s,t_{0})\in  T$ and all $(x,v)\in  Y$.
\end{definition}

\begin{example}\rm\label{ex_ues}
We consider $\mathcal{C}=\mathcal{C}(\mathbb{R}_{+}, \mathbb{R})$ the set of all continuous
functions $x:\mathbb{R}_{+}\rightarrow\mathbb{R}$, with the topology of uniform convergence on
bounded sets. Let $f:\mathbb{R}_{+}\rightarrow \mathbb{R}$ be a decreasing function on $[0,\infty)$.

Let $X$ be the closure in $\mathcal{C}$ of the set of functions
$f_{\theta}, \ \theta\in\mathbb{R}_{+}$ given by
\[
f_{\theta}(\tau)=f(\theta+\tau), \
\forall \tau\in\mathbb{R}_{+}.
\]
The mapping
\begin{equation*}
\varphi: T\times  X\rightarrow  X, \ \varphi(t,s,x)=x_{t-s}
\end{equation*}
is an evolution semiflow on the metric space $ X$, considered with the metric given in Example \ref{ex_cev}.

We consider the Banach space $ V=\mathbb{R}$. The mapping $\Phi: T\times
X \rightarrow \mathcal{B}( \mathbb{R})$ given by
\[
\Phi(t,s,x)v=ve^{-\mu(t-s)+\int_{s}^{t}x(\tau-s)d\tau}, \ (t,
s)\in  T, \ (x,v)\in
 X\times\mathbb{R},
\]
where $\mu>x(0)$, is an evolution cocycle over the evolution semiflow $\varphi$.

The skew-evolution semiflow $C=(\varphi, \Phi)$ is uniformly
exponentially stable with
\[
N=1 \ \textrm{and} \ \nu=\mu-x(0).
\]
\end{example}

\begin{remark}\rm\label{obs_ues_us_ueg}
A uniformly exponentially stable skew-evolution semiflow is
uniformly stable which further implies that it has uniform
exponential growth.
\end{remark}

There exist uniformly stable skew-evolution semiflows which are
not uniformly exponentially stable.

\begin{example}\rm\label{ex_us_neues}
We consider $ X=\mathbb{R}_{+}$, $ V=\mathbb{R}$ and a
nondecreasing and bounded function
$f:\mathbb{R}_{+}\rightarrow\mathbb{R}_{+}^{*}$.

One can easily verify that the mapping $\varphi: T\times
\mathbb{R}_{+}\rightarrow \mathbb{R}_{+}$ defined by
\[
\varphi(t,s,x)=t-s+x, \ (t,s,x)\in T\times \mathbb{R}_{+}
\]
is an evolution semiflow on $\mathbb{R}_{+}$ and the mapping
$\Phi: T\times \mathbb{R}_{+} \rightarrow \mathcal{B}(\mathbb{R})$
given by
\[
\Phi(t,s,x)=\frac{f(x)}{f(t-s+x)}, \ (t,s,x)\in T\times
\mathbb{R}_{+}
\]
is an evolution cocycle on $\mathbb{R}$ and $C=(\varphi, \Phi)$ is
a uniformly stable skew-evolution semiflow, with $N=1$.

On the other hand, if we suppose that there exist some constants
$N\geq 1$ and $\nu >0$ such that
\[
\frac{f(x)}{f(t-s+x)} \leq Ne^{-\nu (t-s)}, \ \forall (t,s,x)\in
T\times \mathbb{R}_{+}.
\]
Let $s=0$ and we obtain
\[
\frac{e^{\nu t}}{f(t+x)}\leq \frac{N}{f(x)}
\]
which is absurd for $t\rightarrow\infty$. Hence $C$ is not
uniformly exponentially stable.
\end{example}

A characterization of the stability property is given by the next
result.

\begin{proposition}\label{caract_ues}
A skew-evolution semiflow $C=(\varphi, \Phi)$ is uniformly
exponentially stable if and only if there exists a nondecreasing
function defined as $f:[0,\infty )\rightarrow (0,\infty )$ with
the property
$\underset{t\rightarrow\infty}\lim f(t)=\infty$ such that%
\begin{equation}
f(t-s)\left\Vert \Phi(t,t_{0},x)v\right\Vert \leq \left\Vert
\Phi(s,t_{0},x)v\right\Vert
\end{equation}%
for all $(t,s),(s,t_{0})\in T$ and all $(x,v)\in
 Y.$
\end{proposition}

\begin{proof}
\emph{Necessity}. It is obvious if we consider $f(t)=N^{-1}e^{\nu
t}$, where $N$ and $\nu$ are given by Definition \ref{ues}.

\emph{Sufficiency}. From the definition of function $f$, there
exists $\delta>0$ such that $f(\delta )>1.$ We denote
\[
\nu =\frac{\ln f(\delta )}{\delta }>0
\]
Let $(t,s)\in T$. Then there exist $n\in\mathbb{N}$ and
$r\in[0,\delta)$ such that $t-s=n\delta+r$. We obtain
\begin{equation*}
e^{\nu (t-s)}\left\Vert \Phi(t,t_{0},x)v\right\Vert \leq f(\delta
)[f(\delta )]^{n}\left\Vert \Phi(t,t_{0},x)v\right\Vert \leq
\end{equation*}%
\begin{equation*}
\leq f(\delta )[f(\delta )]^{n-1}\left\Vert \Phi(t-\delta
,t_{0},x)v\right\Vert \leq ...\leq f(\delta )\left\Vert
\Phi(t-n\delta ,t_{0},x)v\right\Vert \leq
\end{equation*}%
\begin{equation*}
\leq f(\delta )f(r)\left\Vert \Phi(t-n\delta ,t_{0},x)v\right\Vert
\leq f(\delta )\left\Vert \Phi(s,t_{0},x)v\right\Vert.
\end{equation*}
If we consider $N=f(\delta )>1$ it follows that $C$ is uniformly
exponentially stable.
\end{proof}

\vspace{3mm}

Following result can be considered a sufficient condition for a
skew-evolution semiflow to be uniformly exponentially stable. We
will prove a particular case of a result presented by M. Megan,
A.L. Sasu and B. Sasu in \cite{MeSaSa_MIR} for evolution families
and C. Preda and P. Preda in \cite{PrPr_MIR07} for evolution
processes.

\begin{lemma}\label{caract_ues_jum}
Let $C=(\varphi,\Phi)$ be a skew-evolution semiflow with uniform exponential growth with the
property that there exists $\tau>s+1$ such that
\begin{equation}
\left\Vert \Phi(\tau,s,x)v\right\Vert\leq\frac{1}{2}, \ \forall
(s,x,v)\in\mathbb{R}_{+}\times Y \ \textrm{with} \ \left\Vert
v\right\Vert = 1.
\end{equation}
Then $C$ is uniformly
exponentially stable.
\end{lemma}

\begin{proof}
According to the hypothesis there exists $\delta\geq 2$ such that
following relation holds
\[
\left\Vert \Phi(s+\delta,s,x)v\right\Vert<\frac{1}{2}\left\Vert
v\right\Vert, \ \forall s \geq 0, \ \forall(x,v)\in  Y.
\]
Let us consider $(t,s)\in  T$. Then there exist $n\in\mathbb{N}$
and $r\in [0,\delta)$ such that $t=s+n\delta+r$.

We obtain successively
\[
\left\Vert \Phi(t,s,x)v\right\Vert=\left\Vert
\Phi(t,s+n\delta,\varphi(s+n\delta,s,x))\Phi(s+n\delta,s,x)v\right\Vert
\leq
\]
\[
\leq Me^{\omega r}\left\Vert\Phi(s+n\delta,s,x)v\right\Vert
\leq\frac{1}{2}Me^{\omega\delta}\left\Vert\Phi(s+(n-1)\delta,s,x)v\right\Vert\leq
...\leq
\]
\[
\leq 2^{-n}Me^{\omega\delta}\left\Vert v\right\Vert,
\]
for all $(t,s)\in  T$ and for all $(x,v)\in  Y$, where $M$ and
$\omega$ are given by means of Definition \ref{ueg}.

It follows that
\[
\left\Vert \Phi(t,s,x)v\right\Vert\leq Ne^{-\nu(t-s)}\left\Vert
v\right\Vert, \ \forall (t,s)\in  T, \ \forall (x,v)\in  Y,
\]
where we have denoted
\[
N=Me^{2\omega\delta} \ \textrm{and} \ \nu=\delta^{-1}\ln 2,
\]
which ends the proof.
\end{proof}

\begin{theorem}\label{th_RD_toptare}
Let $C=(\varphi,\Phi)$ be a strongly measurable skew-evolution semiflow with uniform
exponential growth. Following statements are equivalent:

$(i)$ $C$ is uniformly exponentially stable

$(ii)$ there exist a function $R\in\mathcal{R}$ and a constant
$N\geq 1$ such that
\begin{equation}
\int_{t}^{\infty}R\left(\left\Vert \Phi(s,t_{0},x)v\right\Vert)
ds\leq NR(\left\Vert \Phi(t,t_{0},x)v\right\Vert\right),
\end{equation}
for all $(t,t_{0},x,v)\in T\times Y$

$(iii)$ there exist a function $R\in\mathcal{R}$ and a constant
$N\geq 1$ such that
\begin{equation}
\int_{t_{0}}^{\infty}R\left(\left\Vert
\Phi(s,t_{0},x)v\right\Vert) ds\leq NR(\left\Vert
v\right\Vert\right),
\end{equation}
for all $(t_{0},x,v)\in\mathbb{R}_{+}\times Y$.
\end{theorem}

\begin{proof}
$(i)\Rightarrow(ii)$ It is obtained immediately by means of Definition \ref{ues} if we consider $R(t)=t, \ t\geq 0$.

$(ii)\Rightarrow(iii)$ It follows by considering $t=t_{0}$.

$(iii)\Rightarrow(i)$ Let $t\geq t_{0}+1$. For $s\in[t-1,t]$ the
existence of constants $M\geq 1$ and $\omega>0$ is assured by
Definition \ref{ueg} such that
\[
\left\Vert\Phi(t,t_{0},x)v\right\Vert\leq Me^{\omega}\left\Vert
\Phi(s,t_{0},x)v\right\Vert, \ \forall (x,v)\in  Y
\]
and, further
\[
R(\left\Vert\Phi(t,t_{0},x)v\right\Vert)\leq
R(Me^{\omega}\left\Vert \Phi(s,t_{0},x)v\right\Vert), \ \forall
(x,v)\in  Y.
\]
By integrating on $[t-1,t]$ we obtain for all $(x,v)\in
 Y$
\[
R(\left\Vert\Phi(t,t_{0},x)v\right\Vert)\leq
\int^{t}_{t-1}R(Me^{\omega}\left\Vert
\Phi(s,t_{0},x)v\right\Vert)ds\leq
\]
\[
\leq
\int_{t_{0}}^{\infty}R(Me^{\omega}\left\Vert\Phi(s,t_{0},x)v\right\Vert)ds\leq
NR(Me^{\omega}\left\Vert v\right\Vert).
\]
We have proved
\[
R(\left\Vert \Phi(t,t_{0},x)v\right\Vert)\leq
NR(Me^{\omega}\left\Vert v\right\Vert), \ \forall t_{0}\geq 0, \
t\geq t_{0}+1, \ \forall (x,v)\in Y.
\]
By integrating on $[t_{0}+1,t]$ it follows
\[
(t-t_{0}-1)R(\left\Vert \Phi(t,t_{0},x)v\right\Vert)=
\]
\[
=\int_{t_{0}+1}^{t}R(\left\Vert
\Phi(t,s,\varphi(s,t_{0},x))v\right\Vert\left\Vert
\Phi(s,t_{0},x)v\right\Vert)ds\leq
\]
\[
\leq N\int_{t_{0}+1}^{t}R(Me^{\omega}\left\Vert
\Phi(s,t_{0},x)v\right\Vert)ds\leq N^{2}R(Me^{\omega}\left\Vert
v\right\Vert).
\]
Hence we have showed
\[
R(\left\Vert \Phi(t,t_{0},x)v\right\Vert)\leq
\frac{N^{2}}{t-t_{0}-1}R(Me^{\omega}),
\]
for all $t>t_{0}+1, \ t_{0}\geq 0$, all $x\in X$ and all $v\in V \
\textrm{with} \ \left\Vert v\right\Vert=1$.

Let $t_{1}>t_{0}+1$ such that
\[
\frac{N^{2}}{t_{1}-t_{0}-1}R(Me^{\omega})<R
\left(\frac{1}{2}\right)
\]
and, further
\[
\left\Vert \Phi(t_{1},t_{0},x)v\right\Vert\leq\frac{1}{2},
\]
for all $x\in X$ and all $v\in V \ \textrm{with} \ \left\Vert
v\right\Vert=1$. Hence, according to Lemma \ref{caract_ues_jum}, $C$ is uniformly
exponentially stable.
\end{proof}

\vspace{3mm}

If we take $R(t)=t^{p}, \ t\geq 0,\ p>1$ we obtain following
result.

\begin{corollary}\label{cor_D_p}
A strongly measurable skew-evolution semiflow $C$ with uniform exponential growth is
uniformly exponentially stable if and only if there exist some
constants $p>1$ and $\widetilde{N}\geq 1$ such that
\[
\int_{t_{0}}^{\infty}\left\Vert \Phi(s,t_{0},x)v\right\Vert^{p}
ds\leq \widetilde{N}^{p}\left\Vert v\right\Vert^{p}, \ \forall
(t,t_{0},x,v)\in T\times Y.
\]
\end{corollary}

\begin{remark}\rm
$(i)$ For $R(t)=t, \ t\geq 0$, from Theorem \ref{th_RD_toptare} is
obtained Theorem 11, proved by Datko in 1972 in \cite{Da_JMA} for
the uniform exponential stability in the strong topology of
$\mathcal{B}( V)$.

$(ii)$ The case of an evolution operator is presented
in 1986 by S. Rolewicz in \cite{Ro_JMAA}, where it is considered a
family of uniformly bounded operators and a nondecreasing function
$N$ vanishing in $0$ but with no convexity restrictions.
\end{remark}

\begin{theorem}\label{th_RD_topunif}
Let $C=(\varphi,\Phi)$ be a skew-evolution semiflow with uniform
exponential growth and the property that for all
$(t_{0},x,v)\in\mathbb{R}_{+}\times Y$ the mappings
$s\mapsto\left\Vert\Phi(s,t_{0},x)v\right\Vert$ and
$s\mapsto\left\Vert\Phi(s,t_{0},x)\right\Vert$ are measurable on
$[t_{0},\infty)$.

Following statements are equivalent:

$(i)$  $C$ is uniformly exponentially stable

$(ii)$ there exist a function $R\in\mathcal{R}$ and a constant
$N\geq 1$ such that
\begin{equation}
\int^{\infty}_{t} R\left(\left\Vert\Phi(s,t_{0},x)\right\Vert
\right)ds\leq N R\left(\left\Vert\Phi(t,t_{0},x)\right\Vert\right)
\end{equation}
for all $(t,t_{0},x)\in T\times X$

$(iii)$ there exists a function $R\in\mathcal{R}$ such that
\begin{equation}
\int^{\infty}_{t_{0}}
R\left(\left\Vert\Phi(s,t_{0},x)\right\Vert\right)ds<\infty
\end{equation}
for all $(t_{0},x)\in\mathbb{R}_{+}\times X$.
\end{theorem}

\begin{proof}
The implication $(i)\Rightarrow(ii)$ follows
easily for $R(t)=t, \ t\geq 0$ and according to Definition \ref{ues}, respectively $(ii)\Rightarrow(iii)$ is obtained for $t=t_{0}$.

$(iii)\Rightarrow(i)$ We have
\[
\int^{\infty}_{t_{0}}
R\left(\left\Vert\Phi(s,t_{0},x)v\right\Vert\right)ds\leq
\int^{\infty}_{t_{0}}
R\left(\left\Vert\Phi(s,t_{0},x)\right\Vert\left\Vert
v\right\Vert\right)ds\leq
\]
\[
\leq\int^{\infty}_{t_{0}}
R\left(\left\Vert\Phi(s,t_{0},x)\right\Vert\right)ds<\infty
\]
for all $(t_{0},x)\in\mathbb{R}_{+}\times X$ and $v\in V$ with
$\left\Vert v\right\Vert\leq 1$. It follows according to Theorem \ref{th_RD_toptare} that $C$ is
uniformly exponentially stable.
\end{proof}

\begin{theorem}\label{th_RB_toptare}
Let $C=(\varphi,\Phi)$ be a $*$-strongly measurable skew-evolution semiflow. Following
statements are equivalent:

$(i)$ $C$ is uniformly exponentially stable

$(ii)$ $C$ is uniformly stable and there exist a function
$R\in\mathcal{R}$ and a constant $\widetilde{N}\geq 1$ such that
\begin{equation}
\int_{t_{0}}^{t}R\left(\left\Vert
\Phi(t,s,\varphi(s,t_{0},x))^{*}v^{*}\right\Vert\right) ds\leq
R\left( \widetilde{N}\left\Vert v^{*}\right\Vert\right),
\end{equation}
for all $(t,t_{0},x,v^{*})\in T\times X\times\mathcal{V^{*}}$

$(iii)$ $C$ has uniform exponential growth and there exist a
function $R\in\mathcal{R}$ and a constant $\widetilde{N}\geq 1$
such that
\begin{equation}
\int_{t_{0}}^{t}R\left(\left\Vert
\Phi(t,s,\varphi(s,t_{0},x))^{*}v^{*}\right\Vert\right) ds\leq
R\left( \widetilde{N}\left\Vert v^{*}\right\Vert\right),
\end{equation}
for all $(t,t_{0},x,v^{*})\in T\times X\times\mathcal{V^{*}}$.
\end{theorem}

\begin{proof}
The implications $(i)\Rightarrow(ii)\Rightarrow(iii)$ are obtained
by considering $R(t)=t, \ t\geq 0$, and according to Remark
\ref{obs_ues_us_ueg}. We obtain also, as $C$ is uniformly exponentially stable
\[
\int_{t_{0}}^{t}\left\Vert
\Phi(t,s,\varphi(s,t_{0},x))^{*}v^{*}\right\Vert ds\leq
N\int_{t_{0}}^{t}e^{-\nu(t-s)}\left\Vert v^{*}\right\Vert ds\leq
\widetilde{N}\left\Vert v^{*}\right\Vert,
\]
for all $(t,t_{0},x,v^{*})\in T\times X\times\mathcal{V^{*}}$,
where we have denoted
\[
\widetilde{N}=\frac{N}{\nu}
\]
and where the constants $N$ and $\nu$ are given as in Definition
\ref{ues}.

$(iii)\Rightarrow(i)$ As a first step we prove that $C$ is
uniformly stable. We consider $t\geq t_{0}+1$ and $s\in[t_{0},t_{0}+1)$. Then
\[
R\left(\frac{1}{M}e^{-\omega(s-t_{0})}|\langle
v^{*},\Phi(t,t_{0},x)v\rangle|\right)\leq
\]
\[
\leq
\int^{t_{0}+1}_{t_{0}}R\left(\frac{\left\Vert\Phi(t,\tau,\varphi(\tau,t_{0},x))^{*}v^{*}%
\right\Vert\left\Vert\Phi(\tau,t_{0},x)v\right\Vert}{Me^{\omega(s-t_{0})}}\right)d\tau\leq
\]
\[
\leq\int^{t_{0}+1}_{t_{0}}R\left(\left\Vert\Phi(t,\tau,\varphi(\tau,t_{0},x))^{*}v^{*}%
\right\Vert \right)d\tau\leq R\left( \widetilde{N}\left\Vert
v^{*}\right\Vert\right),
\]
where the constants $M$ and $\omega$ are given by Definition
\ref{ueg}.

Further we obtain
\[
\left\Vert\Phi(t,t_{0},x)v \right\Vert\leq
Me^{\omega(s-t_{0})}\left\Vert v\right\Vert\leq
Me^{\omega}\left\Vert v\right\Vert, \ \forall t\geq t_{0}+1, \
\forall x\in X.
\]
For $t\in [t_{0},t_{0}+1)$ we have
\[
\left\Vert\Phi(t,t_{0},x) \right\Vert\leq Me^{\omega(t-t_{0})}\leq
Me^{\omega},  \ \forall x\in X.
\]
Hence, $C$ is uniformly stable.

As a second step we prove that $C$ is uniformly exponentially
stable. We remark that for all $t_{0}\geq 0$ there exists
$t>t_{0}$ such that
\[
\frac{N}{t-t_{0}}<R\left(\frac{1}{2M}\right).
\]
For $t>t_{0}$ it follows that
\[
(t-t_{0})R\left(\frac{|\langle
v^{*},\Phi(t,t_{0},x)v\rangle|}{M}\right)\leq
\]
\[
\leq \int^{t}_{t_{0}}R\left(\frac{\left\Vert\Phi(t,s,\varphi(s,t_{0},x))^{*}v^{*}%
\right\Vert\left\Vert\Phi(s,t_{0},x)v \right\Vert}{M}\right)ds\leq
\]
\[
\leq\int^{t}_{t_{0}}R\left(\left\Vert\Phi(t,s,\varphi(s,t_{0},x))^{*}v^{*}\right\Vert\right)ds\leq
N.
\]
Hence, for all $t_{0}\geq 0$ there exists $t>t_{0}$ such that
\[
\frac{|\langle v^{*},\Phi(t,t_{0},x)v\rangle|}{M}\leq
\frac{1}{2M}.
\]
As $\left\Vert v^{*}\right\Vert\leq 1$ we obtain that
\[
\left\Vert\Phi(t,t_{0},x)\right\Vert<\frac{1}{2}, \ \forall x\in
 X.
\]
According to Lemma \ref{caract_ues_jum}, $C$ is uniformly
exponentially stable, which ends the proof.
\end{proof}

\vspace{3mm}

For $R(t)=t^{p}, \ t\geq 0,\ p>1$ we obtain following result.

\begin{corollary}\label{cor_B_p}
A $*$-strongly measurable skew-evolution semiflow $C=(\varphi, \Phi)$ with uniform
exponential growth is uniformly exponentially stable if and only
if there exist some constants $p>1$ and $\widetilde{N}\geq 1$ such
that
\begin{equation}
\int_{t_{0}}^{t}\left\Vert
\Phi(t,s,\varphi(s,t_{0},x))^{*}v^{*}\right\Vert^{p} ds\leq
\widetilde{N}\left\Vert v^{*}\right\Vert^{p},
\end{equation}
for all $(t,t_{0},x,v^{*})\in T\times X\times\mathcal{V^{*}}$.
\end{corollary}

\begin{remark}\rm
For the particular case of exponentially bounded evolution
families of bounded operators acting on Banach spaces similar
results have been obtain by C. Bu\c{s}e, M. Megan, M. Prajea and
P. Preda in \cite{BuMePrPr_IEOT}.
\end{remark}

\begin{theorem}\label{th_RB_topunif}
Let $C=(\varphi,\Phi)$ be a $*$-strongly measurable skew-evolution
semiflow. Following statements are equivalent:

$(i)$ $C$ is uniformly exponentially stable

$(ii)$ $C$ is uniformly stable and there exist a function
$R\in\mathcal{R}$ and a constant $\overline{N}\geq 1$ such that
\begin{equation}
\int^{t}_{t_{0}}R(\left\Vert\Phi(t,s,\varphi(s,t_{0},x))^{*}v^{*}\right\Vert)ds\leq
\overline{N},
\end{equation}
for all $(t,t_{0},x,v^{*})\in T\times X\times V^{*} $ with
$\left\Vert v^{*}\right\Vert\leq 1$.

$(iii)$ $C$ has uniform exponential growth and there exist a
function $R\in\mathcal{R}$ and a constant $\overline{N}\geq 1$
such that
\begin{equation}
\int^{t}_{t_{0}}R(\left\Vert\Phi(t,s,\varphi(s,t_{0},x))^{*}v^{*}\right\Vert)ds\leq
\overline{N},
\end{equation}
for all $(t,t_{0},x,v^{*})\in T\times X\times V^{*} $ with
$\left\Vert v^{*}\right\Vert\leq 1$.
\end{theorem}

\begin{proof}
\emph{Necessity.} Let us consider $R(t)=t$, $t\geq 0$. As $C$ is
uniformly exponentially stable there exist $N\geq 1$ and $\nu>0$
such that
\[
\int_{t_{0}}^{t}\left\Vert\Phi(t,s,\varphi(s,t_{0},x))^{*}v^{*}\right\Vert
ds\leq N\int_{t_{0}}^{t}e^{-\nu(t-s)}\left\Vert v^{*}\right\Vert
ds=
\]
\[
=N\int_{0}^{t-t_{0}}e^{-\nu\tau}\left\Vert v^{*}\right\Vert ds\leq
\overline{N}
\]
for all $t\geq t_{0}+1$ and all $(t,t_{0},x,v^{*})\in T\times
X\times V^{*}$ with $\left\Vert v^{*}\right\Vert\leq 1$, where we
have denoted $\overline{N}=N\nu^{-1}$.

\emph{Sufficiency.} We will first prove that $C$ is uniformly
stable.

Without any loss of generality we can consider that there exists
$K> 0$ such that $R(K)\geq 1$.

Let $N\geq 1$. We consider $t\geq t_{0}+N+1$ and
$s\in[t_{0},t_{0}+N+1)$. As in the proof of Theorem
\ref{th_RB_toptare} we obtain
\[
(N+1)R\left(\frac{1}{M}e^{-\omega(s-t_{0})}|\langle
v^{*},\Phi(t,t_{0},x)v\rangle|\right)\leq
\]
\[
\leq\int^{t_{0}+N+1}_{t_{0}}R\left(\left\Vert\Phi(t,\tau,\varphi(\tau,t_{0},x))^{*}v^{*}%
\right\Vert \right)d\tau\leq N<N+1\leq (N+1)R(K),
\]
where the constants $M$ and $\omega$ are given as in Definition
\ref{ueg}. By taking supremum over $\left\Vert v\right\Vert\leq 1$
and $\left\Vert v^{*}\right\Vert\leq 1$ we obtain
\[
\left\Vert\Phi(t,t_{0},x) \right\Vert\leq
MKe^{\omega(s-t_{0})}\leq MKe^{\omega(N+1)}, \ \forall t\geq
t_{0}+N+1, \ \forall x\in X.
\]
For $t\in [t_{0},t_{0}+N+1)$ we have
\[
\left\Vert\Phi(t,t_{0},x) \right\Vert\leq Me^{\omega(t-t_{0})}\leq
Me^{\omega(N+1)},  \ \forall x\in X
\]
and further
\[
\left\Vert\Phi(t,s,x) \right\Vert\leq M(1+K)e^{\omega(N+1)}, \
\forall (t,s,x)\in T\times X,
\]
which shows that $C$ is uniformly stable.

By a similar argumentation as in Theorem \ref{th_RB_toptare} it
follows that $C$ is uniformly exponentially stable.
\end{proof}

\begin{remark}\rm
Similar results were obtained by C. Bu\c{s}e in \cite{Bu_SED} for
functions of type $R\in\mathcal{R}$ supposed to be convex.
\end{remark}

Often it is possible to formulate results on the asymptotic
behaviors of skew-evolution semiflows in a similar way both for
the continuous and the discrete case.

\begin{proposition}\label{caract_ues_discret}
A skew-evolution semiflow $C=(\varphi,\Phi)$ with uniform
exponential growth is uniformly exponentially stable if and only
if there exist some constants $\widetilde{N}\geq 1$ and
$\widetilde{\nu}> 0$ such that
\begin{equation}
\left\Vert\Phi(n,n_{0},x)v\right\Vert\leq
\widetilde{N}e^{-\widetilde{\nu}(n-n_{0})}\left\Vert v\right\Vert,
\ \forall (n,n_{0},x,v)\in\Delta\times Y.
\end{equation}
\end{proposition}

\begin{proof}
\emph{Necessity}. It can be easily verified, according to
Definition \ref{ues}.

\emph{Sufficiency}. As a first step we will consider $t\geq
t_{0}+1$ and we will denote $n=[t]$ and $n_{0}=[t_{0}]$. This
implies that
\[
n\leq t<n+1, \ n_{0}\leq t_{0}<n_{0}+1, \ n_{0}+1\leq n
\]
We obtain successively
\[
\left\Vert\Phi(t,t_{0},x)v\right\Vert=
\]
\[
=\left\Vert\Phi(t,n,\varphi(n,n_{0}+1,x))\Phi(n,n_{0}+1,\varphi(n_{0}+1,t_{0},x))\Phi(n_{0}+1,t_{0},x)v\right\Vert\leq
\]
\[
\leq
Me^{\omega(t-n)}\left\Vert\Phi(n,n_{0}+1,\varphi(n_{0}+1,t_{0},x))\Phi(n_{0}+1,t_{0},x)v\right\Vert\leq
\]
\[
\leq
Me^{\omega}Me^{\omega(n_{0}+1-t_{0})}\left\Vert\Phi(n,n_{0}+1,x)v\right\Vert
\leq
M^{2}\widetilde{N}e^{2(\omega+\widetilde{\nu})}e^{-\widetilde{\nu}(t-t_{0})}\left\Vert
v\right\Vert,
\]
for all $(x,v)\in  Y$, where the existence of the constants $M$
and $\omega$ is assured by Definition \ref{ueg}.

As a second step, for $t\in[t_{0},t_{0}+1)$ we have
\[
\left\Vert\Phi(t,t_{0},x)v\right\Vert\leq
Me^{\omega(t-t_{0})}\left\Vert v\right\Vert\leq
Me^{\omega+\widetilde{\nu}}e^{-\widetilde{\nu}(t-t_{0})}\left\Vert
v\right\Vert,
\]
for all $(x,v)\in  Y$. Hence, we obtain
\[
\left\Vert\Phi(t,t_{0},x)v\right\Vert\leq \left[
M^{2}\widetilde{N}e^{2(\omega+\widetilde{\nu})}+Me^{\omega+\widetilde{\nu}}\right]e^{-\widetilde{\nu}(t-t_{0})}\left\Vert
v\right\Vert,
\]
for all $(t,t_{0},x,v)\in  T\times Y$, which proves the uniform
exponential stability of the skew-evolution semiflow $C$.
\end{proof}

\vspace{3mm}

A similar result as the sufficient condition for stability
presented in Lemma \ref{caract_ues_jum} can be proved in the discrete time setting.

\begin{lemma}\label{caract_ues_jum_discret}
Let $C=(\varphi,\Phi)$ be a skew-evolution semiflow with uniform exponential growth with the
property that there exists $n_{0}\in\mathbb{N}^{*}$ such that
\begin{equation}
\left\Vert \Phi(n+n_{0},n,x)\right\Vert\leq\frac{1}{2},
\end{equation}
for all $n \in\mathbb{N}, \ \forall x\in  X$. Then $C$ is
uniformly exponentially stable.
\end{lemma}

\begin{proof}
Let us consider $(m,n)\in\Delta$. There exist $k\in\mathbb{N}$ and
$r\in{\{0,1,...,n_{0}-1\}}$ such that $m=n+kn_{0}+r$. Hence
\[
\left\Vert \Phi(m,n,x)\right\Vert\leq
\]
\[
\leq\left\Vert
\Phi(n+kn_{0}+r,n+kn_{0},\varphi(n+kn_{0},n,x))\right\Vert\left\Vert
\Phi(n+kn_{0},n,x)\right\Vert\leq
\]
\[
\leq 2^{-k}Me^{\omega(n_{0}-1)}.
\]
Hence, we obtain
\[
\left\Vert \Phi(m,n,x)\right\Vert\leq
Me^{\omega(n_{0}-1)}e^{-n_{0}^{-1}(m-n)\ln 2}, \ \forall
(m,n)\in\Delta, \ \forall x\in X.
\]
According to Proposition \ref{caract_ues_discret} we obtain the
uniform exponential stability of the skew-evolution semiflow $C$.
\end{proof}

\begin{theorem}\label{th_D_discret}
A skew-evolution semiflow $C=(\varphi,\Phi)$ with uniform
exponential growth is uniformly exponentially stable if and only
if there exist a function $R\in\mathcal{R}$ and a constant
$\widetilde{M}\geq 1$ such that
\begin{equation}
\sum_{n=[t_{0}]+1}^{\infty}R\left(\left\Vert\Phi(n,t_{0},x)v\right\Vert\right)\leq
\widetilde{M}R\left(\left\Vert v\right\Vert\right),
\end{equation}
for all $(n_{0},x,v)\in\mathbb{N}\times Y$ and for all $t_{0}\geq
0$.
\end{theorem}

\begin{proof}
\emph{Necessity}. It follows according to Proposition
\ref{caract_ues_discret} if we consider $R(t)=t, \ t\geq 0$.

\emph{Sufficiency}. First, we will consider $t\geq t_{0}+1$,
$t_{0}\geq 0$. We define $n=[t]$ and $n_{0}=[t_{0}]$. It follows
that
\[
R\left(\left\Vert\Phi(t,t_{0},x)v\right\Vert\right)=R\left
(\left\Vert\Phi(t,n,\varphi(n,t_{0},x))\Phi(n,t_{0},x)v\right\Vert\right)\leq
\]
\[
\leq R\left(Me^{\omega(t-n)}\left\Vert\Phi(n,t_{0},x)v\right\Vert\right)
\leq R\left(Me^{\omega}\left\Vert\Phi(n,t_{0},x)v\right\Vert\right),
\]
for all $(x,v)\in  Y$, where $M$ and $\omega$ are given by
Definition \ref{ueg}.

Further we obtain
\[
\int_{t_{0}+1}^{\infty}R\left(\left\Vert\Phi(t,t_{0},x)v\right\Vert\right)dt
\leq
\sum_{n=[t_{0}]+1}^{\infty}R\left(Me^{\omega}\left\Vert\Phi(n,t_{0},x)v\right\Vert\right)
\leq
\]
\[
\leq M\overline{M}e^{\omega}R\left(\left\Vert v\right\Vert\right),
\]
for all $(x,v)\in  Y$. Then, by Theorem \ref{th_RD_toptare}, there
exist $N\geq 1$ and $\nu>0$ such that
\[
\left\Vert\Phi(t,t_{0},x)v\right\Vert\leq
Ne^{-\nu(t-t_{0})}\left\Vert v\right\Vert, \ \forall t\geq
t_{0}+1, \ \forall (x,v)\in  Y.
\]

We will consider now that $t\in[t_{0},t_{0}+1)$. We have
\[
R\left(\left\Vert\Phi(t,t_{0},x)v\right\Vert\right)\leq
R\left(Me^{\omega(t-t_{0})}\left\Vert v\right\Vert\right)\leq
R\left(Me^{\omega}\left\Vert v\right\Vert\right), \ (x,v)\in
 Y.
\]
As a conclusion we obtain the uniform exponential stability of
$C$.
\end{proof}

\vspace{3mm}

If we consider $R(t)=t^{p}, \ t>0, \ p>0$ we obtain following
result.

\begin{corollary}\label{cor_D_p_discret}
A skew-evolution semiflow $C=(\varphi,\Phi)$ with uniform
exponential growth is uniformly exponentially stable if and only
if there exist $p>0$ and a constant $\widetilde{M}\geq 1$ such
that
\begin{equation}
\left(\sum_{k=n_{0}}^{\infty}\left\Vert\Phi(k,n_{0},x)v\right\Vert^{p}\right)^{\frac{1}{p}}\leq
\widetilde{M} \left\Vert v\right\Vert ,
\end{equation}
for all $(n,x,v)\in\mathbb{N}\times Y$.
\end{corollary}

Another characterization for stability in the discrete time case
is the following result.

\begin{theorem}\label{th_B_discret}
A skew-evolution semiflow $C=(\varphi,\Phi)$ with uniform
exponential growth is uniformly exponentially stable if and only
if there exist a function $R\in\mathcal{R}$ and a constant
$\widetilde{M}\geq 1$ such that
\begin{equation}
\sum_{k=n_{0}}^{\infty}R\left(\left\Vert\Phi(n,k,x)\right\Vert\right)\leq
\widetilde{M},
\end{equation}
for all $(n,n_{0},x)\in\Delta\times X$.
\end{theorem}

\begin{proof}
\emph{Necessity}. It is obtained from the definition of uniform
exponential stability if we consider $R(t)=t,\ t\geq 0$ and
Proposition \ref{caract_ues_discret}.

\emph{Sufficiency}. According to the hypothesis, for all $(k,n_{0})\in\Delta$ we obtain
\[
R\left(\left\Vert\Phi(k,n_{0},x)\right\Vert\right)\leq
\widetilde{M}, \ \forall x\in X.
\]
Further
\[
\sum_{k=n_{0}}^{n}R\left(\left\Vert\Phi(n,n_{0},x)\right\Vert\right)\leq
\sum_{k=n_{0}}^{n}R\left(\left\Vert\Phi(n,k,x)\right\Vert\left\Vert\Phi(k,n_{0},x)\right\Vert\right)\leq
\]
\[
\leq\widetilde{M}
\sum_{k=n_{0}}^{n}R\left(\left\Vert\Phi(n,k,x)\right\Vert\right)\leq
\widetilde{M}^{2}
\]
which implies
\[
R\left(\left\Vert\Phi(n,n_{0},x)\right\Vert\right)
\leq\frac{\widetilde{M}^{2}}{n-n_{0}+1}, \ \forall
(n,n_{0},x)\in\Delta\times X.
\]
There exists $n_{1}>n_{0}+1$ such that
\[
\frac{\widetilde{M}^{2}}{n_{1}-n_{0}+1}\leq
R\left(\frac{1}{2}\right).
\]
We obtain that
\[
\left\Vert\Phi(n_{1},n_{0},x)\right\Vert\leq\frac{1}{2}, \
\textrm{for} \ n_{1}>n_{0}+1 \ \textrm{and} \ x\in X.
\]
According to Lemma \ref{caract_ues_jum_discret}, the conclusion
is obtained.
\end{proof}

\vspace{3mm}

For $R(t)=t^{p}, \ t\geq 0, \ p>0$ following result is obtained.

\begin{corollary}\label{cor_B_p_discret}
A skew-evolution semiflow $C=(\varphi,\Phi)$ with uniform
exponential growth is uniformly exponentially stable if and only
if there exist $p>0$ and a constant $\widetilde{M}\geq 1$ such
that
\[
\left(\sum_{k=n_{0}}^{n}
\left\Vert\Phi(n,k,x)v\right\Vert^{p}\right)^{\frac{1}{p}}\leq
\widetilde{M}\left\Vert v\right\Vert, \ \forall
(n,n_{0},x,v)\in\Delta\times Y.
\]
\end{corollary}

\begin{remark}\rm
In \cite{PrPr_MIR07}, C. Preda and P. Preda have obtained results
that characterize the stability of $\mathcal{C}_{0}$-semigroups
and of evolution processes in continuous and discrete time.
\end{remark}

\section{Exponential stability for skew-evolution semiflows}

In what follows the definitions and characterizations of
exponential growth and exponential stability is obtained in a more
general case, the nonuniform setting. Let $C=(\varphi,\Phi)$ be a
skew-evolution semiflow on $Y= X\times  V$.

\begin{definition}\rm\label{def_neg}
The skew-evolution semiflow $C=(\varphi, \Phi)$ has
\emph{exponential growth} if there exist a couple of applications
$M,\omega:\mathbb{R}_{+}\rightarrow\mathbb{R}_{+}^{*}$ such that
\begin{equation}
\left\Vert \Phi(t,t_{0},x)v\right\Vert \leq M(s)e^{\omega(s)
(t-s)}\left\Vert \Phi(s,t_{0},x)v\right\Vert,
\end{equation}%
for all $(t,s),(s,t_{0})\in  T$ and all $(x,v)\in
 Y$.
\end{definition}

\begin{remark}\rm
In some cases in it more interesting to consider $\omega\equiv c$,
$c$ being a constant, without any loss of generality,
\end{remark}

\begin{remark}\rm\label{cocevol_shift}
Let us define for $\alpha\in\mathbb{R}$ the mapping $C_{\alpha}:
T\times Y\rightarrow  Y$ given by
$C_{\alpha}(t,s,x,v)=(\varphi(t,s,x),\Phi_{\alpha}(t,s,x)v)$,
where $\varphi$ is an evolution semiflow on $ X$ and
\begin{equation}
\Phi_{\alpha}(t,t_{0},x)=e^{-\alpha(t-t_{0})}\Phi(t,t_{0},x), \
(t,t_{0},x)\in T\times X,
\end{equation}
which verifies the conditions of Definition \ref{d_cociclu_2}. It
follows that $C_{\alpha}$ is a skew-evolution semiflow.

Sometimes it is useful to characterize the asymptotic properties
of $C$ by means of those of $C_{\alpha}$

As next relations hold
\[
\left\Vert\Phi_{\alpha}(t,t_{0},x)v\right\Vert=e^{-\alpha(t-t_{0})}
\left\Vert\Phi(t,t_{0},x)v\right\Vert\leq
\]
\[
\leq M(t_{0})e^{[-\alpha+\omega(t_{0})](t-t_{0})}\left\Vert
v\right\Vert
\]
for all $(t,t_{0},x,v)\in T\times Y$, where $M$ and $\omega$ ar
given as in Definition \ref{def_neg}, if follows that $C_{\alpha}$
has also exponential growth.
\end{remark}

\begin{definition}\rm\label{s}
A skew-evolution semiflow $C=(\varphi, \Phi)$ is called
\emph{stable} if there exists a function
$N:\mathbb{R}_{+}\rightarrow\mathbb{R}_{+}^{*}$ such that
\begin{equation}
\left\Vert \Phi(t,t_{0},x)v\right\Vert \leq N(s)\left\Vert
\Phi(s,t_{0},x)v\right\Vert
\end{equation}%
\noindent for all $(t,s),(s,t_{0})\in  T$ and all $(x,v)\in  Y$
\end{definition}

\begin{definition}\rm\label{es}
A skew-evolution semiflow $C=(\varphi, \Phi)$ is said to be
\emph{exponentially stable} if there exist a constant $\nu
>0$ and a mapping
$N:\mathbb{R}_{+}\rightarrow\mathbb{R}_{+}^{*}$ such that
\begin{equation}
\left\Vert \Phi(t,t_{0},x)v\right\Vert \leq N(s)e^{-\nu
(t-s)}\left\Vert \Phi(s,t_{0},x)v\right\Vert,
\end{equation}%
for all $(t,s),(s,t_{0})\in  T$ and all $(x,v)\in
 Y$.
\end{definition}

\begin{example}\rm\label{ex_nues2}
Let $ X=\mathbb{R}_{+}$ and $ V=\mathbb{R}$.

We consider the mapping $\Phi_{f}:
T\times\mathbb{R}_{+}\rightarrow\mathcal{B}(\mathcal{\mathbb{R}})$
given by
\[
\Phi(t,s,x)v=v\exp\left(t\sin t-s\sin s-2t+2s\right), \
(t,s,x,v)\in T\times\mathbb{R}_{+}\times\mathbb{R}.
\]
We obtain that $C=(\varphi,\Phi)$ is a skew-evolution semiflow on
$ Y=\mathbb{R}_{+}\times \mathbb{R}$ over every evolution semiflow
$\varphi: T\times\mathbb{R}_{+}\rightarrow\mathcal{\mathbb{R}}$.

As
\[
t\sin t-s\sin s-2t+2s\leq -t+3s, \ \forall (t,s)\in T
\]
it follows that
\[
\left | \Phi(t,s,x)v\right |\leq e^{2s}e^{-(t-s)}|v|,\ \forall
(t,s,x,v)\in T\times Y,
\]
and, further, that $C$ is exponentially stable.
\end{example}

In what follows, an example  of exponentially stable
skew-evolution semiflows which are not uniformly exponentially
stable will be emphasized.

\begin{example}\rm\label{ex_nues1}
Let us consider $ X=\mathbb{R}_{+}$ and $ V=\mathbb{R}$.

We define the continuous function
$f:\mathbb{R}_{+}\rightarrow[1,\infty)$ given by
\[
f(n)=e^{2n} \ \textrm{and} \
f\left(n+\frac{1}{e^{n^{2}}}\right)=1.
\]
Let us consider the mapping $\Phi_{f}:
T\times\mathbb{R}_{+}\rightarrow\mathcal{B}(\mathcal{\mathbb{R}})$
defined by
\[
\Phi_{f}(t,s,x)v=\frac{f(s)}{f(t)}e^{-(t-s)}v, \ (t,s,x,v)\in
T\times\mathbb{R}_{+}\times\mathbb{R}.
\]
Then $C_{f}=(\varphi,\Phi_{f})$ is a skew-evolution semiflow on $
Y=\mathbb{R}_{+}\times \mathbb{R}$ over any evolution semiflow
$\varphi$ pe $\mathbb{R}_{+}$.

As
\[
\left | \Phi_{f}(t,s,x)v\right |\leq f(s)e^{-(t-s)}|v|, \ \forall
(t,s,x,v)\in T\times Y,
\]
we obtain that $C_{f}$ is exponentially stable.

On the other hand, as
\[
\Phi_{f}(n+\frac{1}{e^{n^{2}}},n,x)=\exp\left(2n-\frac{1}{e^{n^{2}}}\right)\rightarrow\infty
\ \textrm{for} \ n\rightarrow\infty,
\]
it follows that $C_{f}$ is not uniformly exponentially stable.
\end{example}

The next result is a characterization for the property of
exponential stability of skew-evolution semiflows.

\begin{proposition}\label{caract_nues}
A skew-evolution semiflow $C=(\varphi, \Phi)$ is exponentially
stable if and only if there exist a mapping $M:[0,\infty
)\rightarrow (0,\infty )$ and a decreasing function $g:[0,\infty
)\rightarrow (0,\infty )$ with the property
$\underset{t\rightarrow\infty}\lim g(t)=0$ such that
\begin{equation}
\left\Vert \Phi(t,t_{0},x)v\right\Vert \leq M(s)g(t-s)\left\Vert
\Phi(s,t_{0},x)v\right\Vert
\end{equation}%
for all $(t,s,),\ (s,t_{0})\in T$ and all $(x,v)\in
 Y.$
\end{proposition}

\begin{proof}
\emph{Necessity}. It can be easily proven if we consider
\[
M(u)=N(u),\ u\geq 0 \ \textrm{and} \ g(v)=e^{-\nu v}, \ v\geq 0,
\]
where $N$ and $\nu$ are given as in Definition \ref{es}.

\emph{Sufficiency}. As $\underset{t\rightarrow\infty}\lim g(t)=0$,
there exists $\delta\in (0,1)$ such that $g(\delta)<1$. Let
$(t,s)\in T$. Then there exist $n\in\mathbb{N}$ and
$r\in[0,\delta)$ such that $t=s+n\delta+r$. We obtain successively
\[
\left\Vert \Phi(t,t_{0},x)v\right\Vert=\left\Vert
\Phi(s+n\delta+r,t_{0},x)v\right\Vert\leq
\]
\[
\leq g(r) M(s+n\delta)\left\Vert
\Phi(s+n\delta,t_{0},x)v\right\Vert\leq g(0) M(s+n\delta)\left\Vert
\Phi(s+n\delta,t_{0},x)v\right\Vert\leq
\]

\[
\leq g(0)g(\delta)M(s+(n-1)\delta+r)\left\Vert
\Phi(s+(n-1)\delta,t_{0},x)v\right\Vert\leq...\leq
\]
\[
\leq g(0)[g(\delta)]^{n}M(s+(n-1)\delta)...M(s)\left\Vert
\Phi(s+r,t_{0},x)v\right\Vert,
\]
for all $(t,s),\ (s,t_{0})\in T$ and all $(x,v)\in
 Y$.

If we define
\[
N(u)=\frac{g(0)}{g(\delta)}M(u+(n-1)\delta+r)...M(u+r)M(u) \
\textrm{and} \ \nu=\frac{1}{\delta}\ln\frac{1}{g(\delta)}
\]
we obtain that $C$ is exponentially stable.
\end{proof}

\begin{theorem}\label{th_RD_toptare_neunif}
Let $C=(\varphi,\Phi)$ be a skew-evolution semiflow with
exponential growth. Following statements are equivalent:

$(i)$ $C$ is exponentially stable

$(ii)$ there exist a function $R\in\mathcal{R}$, a function
$\widetilde{N}:[0,\infty )\rightarrow (0,\infty )$ and a constant
$\alpha >0$ such that
\begin{equation}
\int_{t}^{\infty}R\left(e^{\alpha(s-t)}\left\Vert
\Phi(s,t_{0},x)v\right\Vert\right) ds\leq
\widetilde{N}(t)R\left(\left\Vert
\Phi(t,t_{0},x)v\right\Vert\right),
\end{equation}
for all $(t,t_{0},x,v)\in T\times Y$

$(iii)$ there exist a function $R\in\mathcal{R}$, a function
$\widetilde{N}:[0,\infty )\rightarrow (0,\infty )$ and a constant
$\alpha >0$ such that
\begin{equation}
\int_{t_{0}}^{\infty}R\left(e^{\alpha(s-t_{0})}\left\Vert
\Phi(s,t_{0},x)v\right\Vert\right) ds\leq
\widetilde{N}(t_{0})R\left(\left\Vert v\right\Vert\right),
\end{equation}
for all $(t_{0},x,v)\in\mathbb{R}_{+}\times Y$.
\end{theorem}

\begin{proof}
$(i)\Rightarrow(ii)$ We consider $R(t)=t, \ t\geq 0$ and let us
define
\[
\alpha=\frac{\nu}{2}>0
\]
$\nu>0$ being given by the fact that $C$ is exponentially stable.

We obtain successively for all $(t,t_{0})\in T$ and $(x,v)\in Y$
\[
\int^{\infty}_{t}e^{\alpha(\tau-t)}\left\Vert
\Phi(\tau,t_{0},x)v\right\Vert
d\tau\leq\int^{\infty}_{t}N(t)e^{\frac{\nu}{2}(\tau-t)}e^{-\nu(\tau-t)}\left\Vert
\Phi(t,t_{0},x)v\right\Vert d\tau=
\]
\[
=\widetilde{N}(t)\left\Vert \Phi(t,t_{0},x)v\right\Vert,
\]
where we have denote
\[
\widetilde{N}(t)=\frac{N(t)}{\alpha}, \ t\geq 0,
\]
the existence of function $N:[0,\infty )\rightarrow (0,\infty )$
being given by Definition \ref{es}.

$(ii)\Rightarrow(iii)$ It is obtained by considering $t=t_{0}$.

$(iii)\Rightarrow(i)$ First, we will prove that the relation in
the hypothesis implies the fact that $C$ is stable. By a similar
proof as in the Sufficiency of Theorem \ref{th_RD_toptare}, we
obtain that
\[
\left\Vert \Phi(t,t_{0},x)v\right\Vert\leq
M(t_{0})N(t_{0})Ke^{\omega(t_{0})}, \ \forall t\geq t_{0}+1, \
\forall (x,v)\in  Y,
\]
where functions $M$ and $\omega$ are given as in Definition
\ref{def_neg} and $K$ is chosen such that $R(K)\geq 1$.

Equally, as $C$ has exponential growth, we have
\[
\left\Vert \Phi(t,t_{0},x)v\right\Vert\leq
M(t_{0})e^{\omega(t_{0})}, \ \forall t\in [t_{0},t_{0}+1), \
\forall (x,v)\in  Y.
\]
Hence we obtain
\[
\left\Vert \Phi(t,t_{0},x)v\right\Vert\leq
M(t_{0})[1+KN(t_{0})]e^{\omega(t_{0})}, \ \forall (t,t_{0})\in T,
\ \forall (x,v)\in  Y,
\]
which proves the stability of $C$.

As a second step we will prove the exponential stability of $C$.

According to the relation in the hypothesis, we obtain that
$C_{-\alpha}$ is stable, where
$C_{-\alpha}=(\varphi,\Phi_{-\alpha})$ is given as in Remark
\ref{cocevol_shift}, which assures the existence of a function
$\widetilde{M}:\mathbb{R}_{+}\rightarrow\mathbb{R}_{+}^{*}$ such
that
\[
e^{\alpha(t-t_{0})}\left\Vert\Phi(t,t_{0},x)v\right\Vert\leq\widetilde{M}(t_{0})\left\Vert
v\right\Vert, \ \forall (t,t_{0},x,v)\in T\times Y,
\]
which proves the exponential stability of the skew-evolution
semiflow $C$ and ends the proof.
\end{proof}

\begin{remark}\rm
For $R(t)=t,\ t\geq 0$ we obtain a generalization in the
nonuniform setting of the classic result proved in Theorem $11$ of
\cite{Da_JMA}.
\end{remark}

\begin{theorem}\label{th_RD_topunif_neunif}
Let $C=(\varphi,\Phi)$ be a skew-evolution semiflow with
exponential growth. Following statements are equivalent:

$(i)$ $C$ is exponentially stable

$(ii)$ there exist a function $R\in\mathcal{R}$, a function
$\widetilde{N}:[0,\infty )\rightarrow (0,\infty )$ and a constant
$\alpha >0$ such that
\begin{equation}
\int_{t}^{\infty}R\left(e^{\alpha(s-t)}\left\Vert
\Phi(s,t_{0},x)\right\Vert\right) ds\leq
\widetilde{N}(t)R\left(\left\Vert
\Phi(t,t_{0},x)\right\Vert\right),
\end{equation}
for all $(t,t_{0},x)\in T\times X$

$(iii)$ there exist a function $R\in\mathcal{R}$ and a constant
$\alpha >0$ such that
\begin{equation}
\int_{t_{0}}^{\infty}R\left(e^{\alpha(s-t_{0})}\left\Vert
\Phi(s,t_{0},x)\right\Vert\right) ds< R(t_{0}),
\end{equation}
for all $(t_{0},x)\in\mathbb{R}_{+}\times X$.
\end{theorem}

\begin{proof}
The implications $(i)\Rightarrow(ii)\Rightarrow(iii)$ follows
easily if we consider $R(t)=t, \ t\geq 0$.

$(iii)\Rightarrow(i)$ As following relation hold
\[
\int_{t_{0}}^{\infty}R\left(e^{\alpha(s-t_{0})}\left\Vert
\Phi(s,t_{0},x)v\right\Vert\right) ds\leq
\]
\[
\leq \left\Vert
v\right\Vert\int_{t_{0}}^{\infty}R\left(e^{\alpha(s-t_{0})}\left\Vert
\Phi(s,t_{0},x)\right\Vert\left\Vert v\right\Vert\right)
ds<R(t_{0})
\]
for all $(t_{0},x,v)\in\mathbb{R}_{+}\times Y$ with $\left\Vert
v\right\Vert\leq 1$, it results by Theorem
\ref{th_RD_toptare_neunif} that $C$ is exponentially stable.
\end{proof}

\begin{theorem}\label{th_RB_neunif}
Let $C=(\varphi,\Phi)$ be a skew-evolution semiflow with
exponential growth. $C$ is exponentially stable if and only if
there exist a function $R\in\mathcal{R}$, a mapping
$N:\mathbb{R}_{+}\rightarrow[1,\infty)$ and a constant $\alpha>0$
such that
\begin{equation}
\int^{t}_{t_{0}}R\left(e^{\alpha(t-s)}\left\Vert\Phi(t,s,\varphi(s,t_{0},x))^{*}v^{*}\right\Vert
\right)ds\leq N(t_{0}),
\end{equation}
for all $(t,t_{0},x,v^{*})\in T\times X\times V^{*}$ with
$\left\Vert v^{*}\right\Vert\leq 1$.
\end{theorem}

\begin{proof} \emph{Necessity}. If we consider $R(t)=t, \ t\geq 0$ and if we define
$\alpha=-\frac{\omega}{2}$ we obtain
\[
\int^{t}_{t_{0}}e^{\alpha(t-s)}\left\Vert\Phi(t,s,\varphi(s,t_{0},x))^{*}v^{*}\right\Vert
ds\leq M(t_{0})\int^{t}_{t_{0}}e^{\frac{\omega}{2}(t-s)}\left\Vert
v^{*}\right\Vert ds\leq N(t_{0}),
\]
for all $(t,t_{0},x,v^{*})\in T\times X\times V^{*}$ with
$\left\Vert v^{*}\right\Vert\leq 1$, where we have defined
\[
N(u)=\frac{2}{\omega}M(u),
\]
function $M$ given by the exponential growth of $C$.

\emph{Sufficiency}. Let us consider $t\geq t_{0}+N(t_{0})+1$. We
have
\[
\left[N(t_{0})+1\right]R\left(\frac{e^{\alpha(t-t_{0})}e^{-\omega(t_{0})(s-t_{0})}}{M(t_{0})}\left|\left\langle
v^{*},\Phi(t,t_{0},x)v\right\rangle \right|\right)
\]
\[
\leq
\int^{t_{0}+N(t_{0})+1}_{t_{0}}R\left(\frac{e^{\alpha(t-s)}\left\Vert\Phi(t,s,\varphi(s,t_{0},x))^{*}v^{*}%
\right\Vert\left\Vert\Phi(s,t_{0},x)v\right\Vert}{M(t_{0})e^{\omega(t_{0})(s-t_{0})}e^{-\alpha(s-t_{0})}}\right)ds\leq
\]
\[
\leq\int^{t_{0}+N(t_{0})+1}_{t_{0}}R\left(e^{\alpha(t-s)}\left\Vert \Phi(t,s,\varphi(s,t_{0},x))^{*}v^{*}%
\right\Vert \right)ds\leq
\]
\[
\leq N(t_{0})< [N(t_{0})+1]R(1),
\]
where functions $M$ and $\omega$ are assured by Definition
\ref{def_neg} and where, without loss of generality, we have
supposed $R(1)>1$. By taking supremum over $\left\Vert
v\right\Vert\leq 1$ and $\left\Vert v^{*}\right\Vert\leq 1$ we
obtain
\[
e^{\alpha(t-t_{0})}\left\Vert\Phi(t,t_{0},x) \right\Vert\leq
M(t_{0})e^{\omega(t_{0})(s-t_{0})}\leq
M(t_{0})e^{\omega(t_{0})[N(t_{0})+1]},
\]
for all $t\geq t_{0}+N(t_{0})+1$ and all $x\in X$.

If we take $t\in [t_{0},t_{0}+N(t_{0})+1)$ we have
\[
\left\Vert\Phi(t,t_{0},x) \right\Vert\leq
M(t_{0})e^{[\omega(t_{0})-\alpha](t-t_{0})}\leq
M(t_{0})e^{[\omega(t_{0})-\alpha][N(t_{0})+1]}
\]
for all $x\in X$.

It follows that $C_{-\alpha}$ is stable, where
$C_{-\alpha}=(\varphi,\Phi_{-\alpha})$ is given as in Remark
\ref{cocevol_shift}. Then there exists a function
$\widetilde{N}:\mathbb{R}_{+}\rightarrow [1,\infty)$ such that
\[
e^{\alpha(t-t_{0})}\left\Vert\Phi(t,t_{0},x)\right\Vert\leq\widetilde{N}(t_{0}),
\ \forall (t,t_{0},x)\in T\times X,
\]
which, as $\alpha>0$ implies the exponential stability of the
skew-evolution semiflow $C$.
\end{proof}

\vspace{3mm}

As in the uniform setting, we can describe asymptotic properties
of skew-evolution semiflows, as stability, in the discrete time
case. Such a characterization for the exponential stability is
emphasized in the next

\begin{proposition}\label{caract_es_discret}
A skew-evolution semiflow $C=(\varphi,\Phi)$ with exponential
growth is exponentially stable if and only if there exist a
constant $\mu> 0$ and a sequence of real numbers $(a_{n})_{n\geq
0}$ with the property $a_{n}\geq 1, \ \forall n\geq 0$ such that
\begin{equation}
\left\Vert\Phi(n,n_{0},x)v\right\Vert\leq
a_{n}e^{-\mu(n-n_{0})}\left\Vert v\right\Vert,
\end{equation}
for all $(n,n_{0},x,v)\in\Delta\times Y$.
\end{proposition}

\begin{proof}
\emph{Necessity}. It can be easily shown.

\emph{Sufficiency}. If $t\geq t_{0}+1$ and $n=[t]$ and
$n_{0}=[t_{0}]$ we have
\[
n\leq t<n+1, \ n_{0}\leq t_{0}<n_{0}+1, \ n_{0}+1\leq n
\]
We obtain
\[
\left\Vert\Phi(t,t_{0},x)v\right\Vert\leq
\]
\[
\leq
M(n)e^{\omega(n)(t-n)}\left\Vert\Phi(n,n_{0}+1,\varphi(n_{0}+1,t_{0},x))\Phi(n_{0}+1,t_{0},x)v\right\Vert\leq
\]
\[
\leq a_{n}M^{2}(n)e^{2[\omega(n)+\mu]}e^{-\mu(t-t_{0})}\left\Vert
v\right\Vert,
\]
for all $(x,v)\in  Y$, where the existence of the functions $M$
and $\omega$ is assured by Definition \ref{def_neg}.

As a second step, for $t\in[t_{0},t_{0}+1)$ we have
\[
\left\Vert\Phi(t,t_{0},x)v\right\Vert\leq
M(t_{0})e^{\omega(t_{0})(t-t_{0})}\left\Vert v\right\Vert\leq
M(t_{0})e^{\omega(t_{0})+\mu} e^{-\mu(t-t_{0})}\left\Vert
v\right\Vert,
\]
for all $(x,v)\in  Y$.

Hence, we obtain the exponential stability of $C$.
\end{proof}

\vspace{3mm}

Some characterization for the stability of skew-evolution
semiflows in the nonuniform setting similar to the continuous case
are presented in discrete time.

\begin{theorem}\label{th_D_discret_neunif}
A skew-evolution semiflow $C=(\varphi,\Phi)$ is exponentially
stable if and only if there exist a function $R\in\mathcal{R}$, a
constant $\rho>0$ and a sequence of real numbers
$(\alpha{n})_{n\geq 0}$ with the property $\alpha_{n}\geq 1, \
\forall n\geq 0$ such that
\begin{equation}
\sum_{k=n}^{m}R\left(e^{\rho(k-n)}\left\Vert
\Phi(k,n,x)v\right\Vert\right)\leq\alpha_{n} R\left(\left\Vert
v\right\Vert\right)
\end{equation}
for all $(m,n,x,v)\in\Delta\times Y$.
\end{theorem}

\begin{proof} \emph{Necessity}. Let us take $R(t)=t, \ t\geq 0$.
Definition \ref{es} provides the constant $\nu> 0$ and the
sequence of real numbers $(a_{n})_{n\geq 0}$ with the property
$a_{n}\geq 1, \ \forall n\geq 0$. We obtain for
$\rho=\frac{\nu}{2}>0$ and according to Proposition
\ref{caract_es_discret}
\[
\sum_{k=n}^{m}e^{\rho(k-n)}\left\Vert\Phi(k,n,x)v\right\Vert \leq
a_{n} \sum_{k=n}^{m}e^{\rho(k-n)}e^{-\nu(k-n)}\left\Vert
\Phi(n,n,x)v\right\Vert =
\]
\[
=a_{n}\left\Vert v\right\Vert
\sum_{k=n}^{m}e^{-\frac{\nu}{2}(k-n)}\leq \alpha_{n}\left\Vert
v\right\Vert, \forall m,n \in \Delta, \forall (x,v)\in Y,
\]
where we have denoted
\[
\alpha_{n}=a_{n}e^{\frac{\nu}{2}}, \ n\in \mathbb{N}.
\]

\emph{Sufficiency}. According to the hypothesis and by a similar
proof as in Theorem \ref{th_D_discret} we have that
$C_{-\rho}=(\varphi,\Phi_{-\rho})$ is stable, where
\[
\Phi_{-\rho}(m,n,x)=e^{\rho(m-n)}\left\Vert\Phi(m,n,x)\right\Vert,
\ (m,n,x)\in\Delta\times X.
\]
Thus there exists a  sequence of real numbers $(a_{n})_{n\geq 0}$
with the property $a_{n}\geq 1, \ \forall n\geq 0$, such that
\[
e^{\rho(m-n)}\left\Vert\Phi(m,n,x)v\right\Vert\leq a_{n}\left\Vert
v\right\Vert, \ \forall (m,n,x,v)\in\Delta\times Y,
\]
which implies the exponential stability of the skew-evolution
semiflow $C$ and ends the proof.
\end{proof}

\begin{theorem}\label{th_B_discret_neunif}
A skew-evolution semiflow $C$ is exponentially stable if and only
if there exist a function $R\in\mathcal{R}$, a constant $\gamma>0$
and a sequence of real numbers $(\beta{n})_{n\geq 0}$ with the
property $\beta_{n}\geq 1, \ \forall n\geq 0$ such that
\begin{equation}
\sum_{k=n}^{m}R\left(e^{\gamma(m-k)}\left\Vert
\Phi(m,k,\varphi(k,n,x))^{*}v^{*}\right\Vert\right)\leq\beta_{n}
R\left(\left\Vert v^{*}\right\Vert\right)
\end{equation}
for all $(m,n,x,v)\in\Delta\times Y$.
\end{theorem}

\begin{proof} \emph{Necessity}. For $R(t)=t, \ t\geq 0$ and
$\gamma=\frac{\nu}{2}>0$ we obtain according to the definition of
stability and Proposition \ref{caract_es_discret}
\[
\sum_{k=n}^{m}e^{\frac{\nu}{2}(m-k)}\left\Vert\Phi(m,k,\varphi(k,n,x))^{*}v^{*}\right\Vert
\leq a_{n}\left\Vert
v^{*}\right\Vert\sum_{k=n}^{m}e^{-\frac{\nu}{2}(m-k)}\leq\beta_{n}
\left\Vert v^{*}\right\Vert
\]
where we have denoted
\[
\beta_{n}=\frac{a_{n}}{1-e^{-\frac{\nu}{2}}}, \ n\in \mathbb{N},
\]
where the constant $\nu> 0$ and the sequence of real numbers
$(a_{n})_{n\geq 0}$ with the property $a_{n}\geq 1, \ \forall
n\geq 0$ are given by Definition \ref{es}.

\emph{Sufficiency}. According to the hypothesis and by a similar
proof as in Theorem \ref{th_B_discret} we have that
$C_{-\gamma}=(\varphi,\Phi_{-\gamma})$ given as in Remark
\ref{cocevol_shift} is stable. It follows as in the proof of
Theorem \ref{th_D_discret_neunif} that the skew-evolution semiflow
$C$ is exponentially stable.
\end{proof}


{\footnotesize

\end{document}